\documentclass[11pt]{amsart}
\usepackage{graphicx}

\newtheorem{thm}{Theorem}
\newtheorem{prop}{Proposition}
\newtheorem{lemma}[prop]{Lemma}
\newtheorem{cor}[prop]{Corollary}
\theoremstyle{definition}
\newtheorem*{defn}{Definition}

\newtheorem*{problem}{Open Problem}

\newcommand{\ignore}[1]{}

\begin{document}

\title{Lower Bounds for Simplicial Covers and Triangulations of Cubes}

\author[Adam Bliss]{Adam Bliss $^*$}
\thanks{$^*$Research partially supported by a Beckman Research Grant at
Harvey Mudd College.}
\address{Department of Mathematics\\ Harvey Mudd College\\ Claremont, CA
91711\\ U.S.A.}
\email{abliss@hmc.edu}
\author[Francis Su]{Francis Edward Su $^{**}$}
\thanks{$^{**}$Research partially supported by NSF Grant DMS-0301129.}
\address{Department of Mathematics\\ Harvey Mudd College\\ Claremont, CA
91711\\ U.S.A.}
\email{su@math.hmc.edu}

\begin{abstract}
We show that the size of a minimal simplicial cover of a polytope $P$
is a lower bound for the size of a minimal triangulation of $P$,
including ones with extra vertices.  We then use this fact to study
minimal triangulations of cubes, and we improve lower bounds for
covers and triangulations in dimensions 4 through at least 12 
(and possibly more dimensions as well).
Important ingredients are an analysis of the number of 
exterior faces that a simplex in the
cube can have of a specified dimension and volume, and a
characterization of corner simplices in terms of their exterior faces.
\end{abstract}

\maketitle


Let $P$ be a convex polytope. A {\em (simplicial) cover} of
$P$ is a collection of simplices such that (i) the vertices of the
simplices are vertices of $P$ and (ii) the union of the
simplices is $P$.  For example, some covers are triangulations
in which
simplices meet face-to-face and have disjoint interiors; although in general, 
cover elements may overlap.

Define the {\em covering number} $C(P)$ to be the 
minimal number of simplices needed for a cover of a polytope $P$.  
Although the covering number is of interest in its own right 
(see \cite{deloera-peterson-su}),
we prove in Theorem \ref{covers-bound}
that the covering number of $P$ also gives a 
lower bound for the size of a minimal triangulation of $P$, 
including triangulations with extra vertices.

We then use the covering number to study the classical problem 
of determining the size of the minimal triangulation when $P$ is a
$d$-dimensional cube (this is sometimes called the {\em simplexity} of the
cube).  Let $I=[0,1]$, and let $I^d$ denote the $d$-cube.
We define the notion of an {\em exterior face} of a cube simplex,
develop a counting function for exterior faces of prescribed dimension
and volume, and establish a recursive bound in 
Theorem \ref{thm:recurrence} that yields an absolute bound for the
cases of interest in Theorem \ref{c-prime=c}.  
This bound can be improved still further
by analyzing separately the corner simplices of the cube; Theorem
\ref{non-corner} characterizes corner simplices in terms of the number
of exterior faces they have.
We use these results to establish new lower bounds for $C(I^d)$, 
and via Theorem \ref{covers-bound}, these yield the best known
bounds for general triangulations in dimensions up through 
at least $12$ (and possibly more dimensions--- see the concluding section).  
For $d=4$ through $d=12$ these numbers are:
$$
16, 60, 252, 1143, 5104, 22616, 98183, 520865, 2.9276 \times 10^6
$$
In particular, our results show that Mara's minimal triangulation
of the 4-cube \cite{mara}, using 16 simplices,
is also a minimal simplicial cover; furthermore,
using extra vertices will not produce any smaller triangulations of the 4-cube.
See Theorem \ref{4-cube-cover}.  
A comparison of our bounds
with bounds for other kinds of decompositions may be found in 
Table \ref{table-bounds}.

This paper is organized as follows.  Section 1 discusses the
relationship between minimal simplicial covers and minimal
triangulations of arbitrary polytopes $P$.  Section 2 gives
background on triangulations of cubes.  We develop some
terminology in Section 3, and in Section 4, we develop
constraints for a linear program for our problem.  This
involves a counting function $F$ which counts the number of exterior
faces that a simplex can have; in order to estimate $F$, we develop in
Sections 5 and 6 some theory regarding the way exterior faces of
simplices relate to simplices in the cube.  Then, in Section 7 we show
how this theory produces a recurrence for the counting function that
can be used to get bounds on $F$, and Section 8 refines the
earlier linear program.  The final section discusses our results and
some open questions.

\section{Minimal Covers bound Triangulations with extra vertices}

When speaking of minimal triangulations of a polytope $P$,
we must be careful to distinguish what kind of decomposition we mean,
since there are many such notions in the literature.
In this paper, a {\em triangulation} of a $d$-polytope $P$ 
will always mean a decomposition of
$P$ into $d$-simplices with disjoint interiors,
such that each pair of simplices intersect in a face common to both
or not at all.  

In the literature,
some authors restrict attention to triangulations
whose vertices are required to be vertices of the polytope.  We 
refer to these as {\em vertex triangulations}.
Others do not always require triangulations to meet face-to-face; such
decompositions are sometimes called {\em dissections}.  They are
covers by simplices that do not have overlapping interiors.  
A {\em vertex dissection} is a dissection in which 
vertices of simplices come from those of the polytope.  

We let $D(P),T(P),D^v(P),T^v(P)$ denote, respectively,
the size of the smallest possible 
dissection, triangulation, vertex dissection,
and vertex triangulation of $P$.  (We note that our notation differs
from that of Smith \cite{smith}, who uses $T$ to include dissections,
and of Hughes-Anderson
\cite{hughes-anderson}, who use $T$ for vertex triangulations.)

The evident inclusions immediately
imply that $D(P) \leq T(P) \leq T^v(P)$ and
$D(P) \leq D^v(P) \leq T^v(P)$.
Also, $C(P) \leq D^v(P)$ since any vertex dissection is a simplicial
cover--- this does not hold for general dissections, since in a 
cover of $P$, vertices of simplices must be vertices of $P$.
What may therefore be somewhat surprising is the following important 
inequality relating covers to general triangulations:

\begin{thm}
\label{covers-bound}
For any convex polytope $P$, the covering number $C(P)$ satisfies:
$$
C(P) \leq T(P).
$$
\end{thm}

Some covers arise as images of simplicial Brouwer self-maps 
of a triangulated polytope and are related to the polytopal Sperner
lemma of DeLoera-Peterson-Su \cite{deloera-peterson-su}. 
Such considerations form the basis for a proof of Theorem \ref{covers-bound}.

\begin{figure}[thpb]
\begin{center}
\includegraphics[height=1.8in]{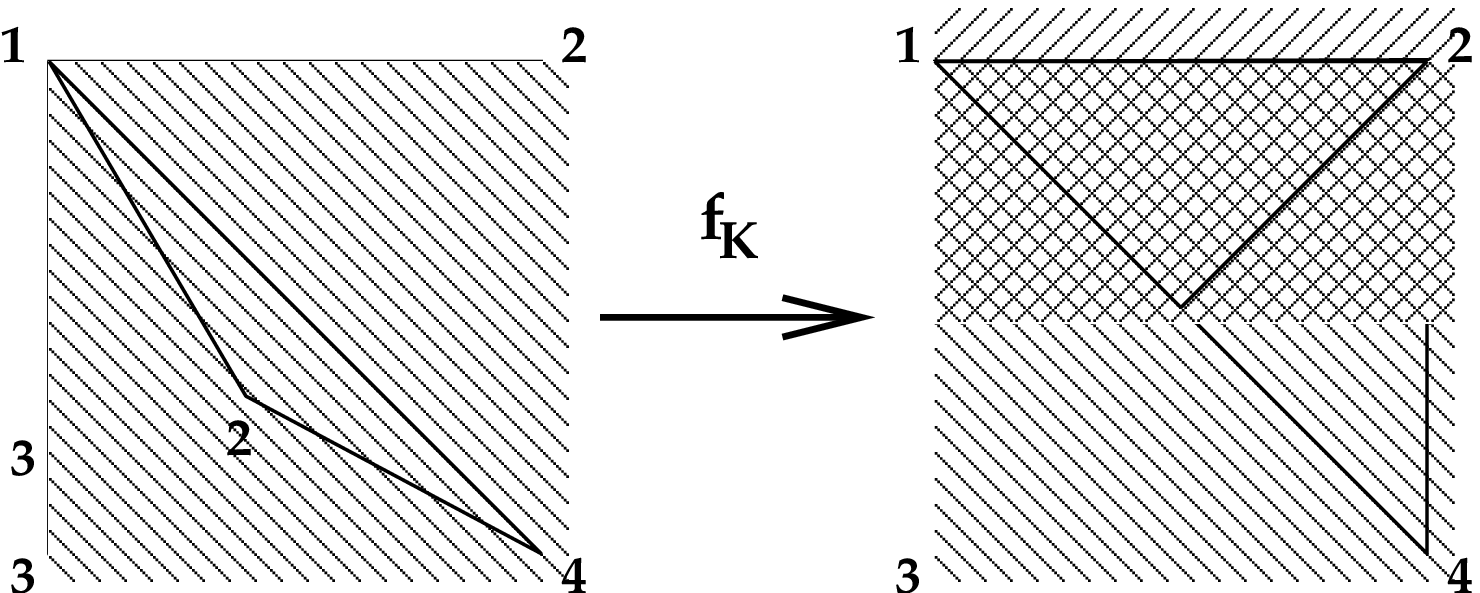}
\end{center}
\caption{The piecewise linear map $f_K$.  
  Two of the simplices have been shaded to show
  how they behave under $f_K$.
Note that the images of simplices in a 
  triangulation can overlap.}
\label{coverbound}
\end{figure}

\begin{proof}
Let $K$ be a triangulation of $P$.  For any such triangulation
construct a piecewise linear (PL) map 
$f_K:P \rightarrow P$ in the following way.  If $v$ is a vertex of $K$,
define $f_K(v)$ to be any vertex of $P$ on the smallest-dimensional
$P$-face that contains $v$.
Then extend this map linearly across each simplex of $K$.  See Figure 
\ref{coverbound}.  

Thus $f_K$ is a PL-map from $P$ to $P$ that takes simplices of $K$ to
simplices formed by vertices of $P$, and these images must be a cover of $P$ 
because this map is a Brouwer map of degree 1
\cite[Prop.~3]{deloera-peterson-su}.  
Thus there are at least $C(P)$ such simplices in the triangulation $K$.
\end{proof}

Thus the covering number is a lower bound for the size of any
triangulation (including ones with extra vertices).  
We make two remarks about the proof.
First, the above assignment of $f_K$ to vertices of $P$ 
is called a {\em Sperner labelling} of the
vertices, in the sense of
\cite{deloera-peterson-su}.  Secondly,
this proof does not work for arbitrary dissections, because if
the simplices do not meet face-to-face, the resulting
PL-map may not be well-defined.  

We may summarize the relationships discussed above as
\begin{equation}
C(P), D(P) \ \leq \ T(P), D^v(P) \ \leq\  T^v(P)
\end{equation}
which signifies a partial ordering in which there is no known relation
between $C$ and $D$ nor between $T$ and $D^v$, but all other
inequalities hold.
Very little is known about the strictness of the inequalities above
for any class of polytopes $P$.
Constructions in \cite{deloera-below} show that in some polytopes, 
$C(P)$ can be strictly smaller than $T^v(P)$, although it is unknown
whether this occurs for cubes.

\section{Triangulations of Cubes}
We now restrict our attention to cubes.
Interest in small triangulations of cubes stems principally from certain
simplicial fixed-point algorithms (e.g. see \cite{todd}) which run faster
when there are fewer simplices.  
The same considerations govern more recent 
{\em fair division procedures} \cite{simmons-su, su} that depend on efficient
triangulations.  

The $d$-cube $[0,1]^d$
has a standard triangulation $T_0$ of size $d!$ 
in which each permutation $(x_{s_1},...,x_{s_n})$ 
of $n$ variables $x_1,..., x_n$ is associated
with the simplex of points for which 
$0 \leq x_{s_1} \leq x_{s_2} \leq x_{s_3} \leq ... \leq x_{s_n} \leq 1$.  
From this description, it is easy to see that the simplices meet only
in common faces and have no common interior points.  This
triangulation is maximal among those that use only vertices of the cube.

However, only for $d=2$ is this triangulation minimal.  
For $d=3$, it is easy to check that the 
$T^v$-minimal triangulation is of size 5, 
formed by four corner simplices at non-adjacent vertices, and one 
fat regular tetrahedron using the other four vertices.  It follows
from results of Smith \cite{smith} that this is minimal for all other
kinds of decompositions as well.
For $d=4$, Mara \cite{mara} produced a triangulation 
using only 16 simplices; 
Cottle \cite{cottle} and Sallee \cite{sallee-82dm} 
showed it was minimal for $T^v(I^4)$ and Hughes \cite{hughes} showed
$D^v(I^4)=16$.  In Theorem \ref{4-cube-cover}, we show that $C(I^4)=16$
and therefore $T(I^4)=16$.  It is unknown whether $D(I^4)=15$ or $16$.
Hughes-Anderson \cite{hughes-anderson} showed that
$T^v(I^5)=67, T^v(I^6)=308, T^v(I^7)=1493$, but for $d \geq 8$ 
there are no exact results for $T^v(I^d)$.  
For $d \geq 5$, there are no exact results known for $C, D, T, D^v$.
We improve the best bounds known for $C$ and $T$ in many small
dimensions.  See Table \ref{table-bounds}.

\begin{table}[t]
\caption{Comparison of best-known 
lower bounds for cubes, and kinds of decompositions
to which they apply.
Equal signs denote the cases for $d\geq 4$ for which bounds are
known to be optimal.
}
\label{table-bounds}
{\footnotesize
\begin{center}
\begin{tabular}{|rrrr|}
\hline
    &            &            & Hughes \cite{hughes},\\
Dimension 
    & Smith \cite{smith}
                 & our bounds& Hughes-Anderson*\cite{hughes-anderson},\\
    &            &           & Cottle**\cite{cottle} and Sallee**\cite{sallee-82dm}\\
    & $C,D,T,D^v,T^v$ &  $C,T,D^v,T^v$ & $D^v,T^v$\\
\hline
3   &  5         &  5       &  5           \\
4   &  15        &  $C,T$=16&  $D^v$=16,$T^v$=16**\\
5   &  48        &  60      &  61, $T^v$=67* \\
6   &  174       &  252     &  270*, $T^v$=308* \\
7   &  681       &  1,143   &  1,175*, $T^v$=1493*\\
8   &  2,863     &  5,104   &  5,522 \\
9   &  12,811    &  22,616  &  26,593 \\
10  &  60,574    &  98,183  &  131,269 \\
11  &  300,956   &  520,865 &  665,272 \\ 
12  &  1,564,340 &  $2.9276 \times 10^6$ & \\
& & & \\
$d$ & ${\mbox{\tiny asymptotic bound} \atop
\frac{6^{\frac{d}{2}}d!}{2(d+1)^{\frac{d+1}{2}}}}$  & & \\
\hline
\end{tabular}
\end{center}
}
\end{table}

To obtain an asymptotic bound for the size of the minimal cover 
of a $d$-cube, one might use the following simple idea: 
the number of simplices in the cover
is bounded below by the volume of the cube
divided by the volume of the largest possible simplex in the
$d$-cube.  This volume is always of the form $V(d)/d!$ for $V(d)$ a
positive integer.
Determining $V(d)$ is a hard problem related to the 
{\em Hadamard determinant problem}, but for small $d$ the values are
known.  See Table \ref{vd-table}.

\begin{table}[h]
\caption {Some values of $V(d)$ from \cite{hudelson-klee-larman}.}
\label{vd-table}
\begin{tabular}{|r|llllllllllllll|}
\hline
$d$ & 0&1&2&3&4&5&6&7&8&9&10&11&12&13\\
\hline
$V(d)$ &1&1&1&2&3&5&9&32&56&144&320&1458&3645&9477\\
\hline
\end{tabular}
\end{table}
The survey \cite{hudelson-klee-larman} also gives some infinite 
families of answers, and an asymptotic upper bound for $V(d)$:
\begin{equation}
  \label{vd-upperbound} 
  V(d) \leq \frac{(d+1)^{\frac{d+1}{2}}}{2^d}
\end{equation}
from which one may obtain $C(d) \geq d! 2^d (d+1)^{-\frac{d+1}{2}}$.
Smith \cite{smith} observed that one can improve this technique by
considering {\em hyperbolic volumes} instead of Euclidean volumes,
and arrived at the improved bound
\begin{equation}
  \label{smith-bound}
  C(d) \ge \frac{6^{\frac{d}{2}}d!}{2(d+1)^{\frac{d+1}{2}}}.
\end{equation}
(Although Smith only spoke of triangulations and non-overlapping
decompositions of the cube, his method actually also
applies to covers as well.)
This asymptotic bound remains the
best asymptotic bound for $C(d)$ for arbitrary $d$,
but our methods improve the explicit bounds that Smith gives 
for low dimensions; see Table~\ref{table-bounds}.

However, his results are not fully comparable to ours, since his
methods also apply to general dissections.  Thus our work gives tighter
lower bounds for covers and triangulations in specific dimensions.

In contrast to Smith, we avoid the use of hyperbolic geometry.  
Instead, we develop a linear
program whose optimal solution is a bound for the size of the minimal cover.
The linear programming approach was initiated by 
Sallee \cite{sallee-82dam}, Hughes \cite{hughes}, and
Hughes-Anderson \cite{hughes-anderson} to study triangulations of the
cube.  However, the results of Hughes and Hughes-Anderson are not
fully comparable to ours, because their methods do not apply to
covers.

On the other hand, Sallee's method does apply to covers (although he 
only spoke of vertex triangulations), and so our Theorem
\ref{covers-bound} shows that his results hold for general triangulations as
well.  As reported by Hughes \cite{hughes}, Sallee's method gives
these lower bounds for $d=3$ though $d=11$:
$
5, 16, 60, 250, 1117, 4680, 21384, 95064, 502289,
$
although the more recent determination of $V(10)$ would improve the 
last two bounds to 95708, 516465.
Thus our bounds agree with his for $d \leq 5$ and dominate
his bounds for $d > 5$.  

Upper bounds for minimal triangulations can be obtained by construction.
Recent work of Orden-Santos \cite{orden-santos} shows that 
the $d$-cube can be triangulated with $O(.816^d d!)$ simplices, so
there remains a large gap between this asymptotic upper bound and the
asympototic lower bound of Smith.  A survey of specific upper bounds in low
dimensions (up through 12) may be found in \cite{smith}.

\section{Exterior faces of simplices in the cube}
Hereafter when we refer to a simplex, we will (unless otherwise
specified) mean a non-degenerate simplex 
spanned by vertices of the unit $d$-cube.

A {\em $j$-face of a $d$-simplex} is the $j$-simplex spanned by
some $j+1$ of the simplex's vertices.  
A {\em $j$-face of a $d$-cube} is the $j$-cube spanned by some 
$2^j$ of the cube's vertices that lie in a $j$-dimensional hyperplane 
on which $d-j$ coordinates agree.  
In both cases the number $d-j$ is called the {\em codimension} of
the face.  
A face of codimension 1 is called a {\em facet}.
Two $j$-faces of a $d$-cube are said to be {\em parallel} if the
$j$-dimensional hyperplanes containing them are parallel.

We say a $j$-face of a $d$-simplex is 
{\em exterior} if it is contained in a $j$-face of the $d$-cube.
The empty set will also be considered an exterior face.
As an example, the diagonal of a facet of the $3$-cube is {\em not} an exterior
$1$-face, because it is not contained in a $1$-face of the $d$-cube.

We can represent a $d$-simplex $\sigma$ in a $d$-cube as a
$(d+1)\times d$ matrix $M$ in which the rows are 
coordinates for the vertices of $\sigma$.  
We shall call $M$ the {\em matrix representation} of $\sigma$.
Let $[1|M]$ denote $(d+1)\times(d+1)$ square matrix formed by
augmenting $M$ by an initial column of ones.  
Then $|\det[1|M]|/d!$ is the volume of $\sigma$.
In particular, since the vertices of $\sigma$ are
chosen from $\{0,1\}^d$, every entry in $M$ is either a zero or a one,
so volumes of $d$-simplices are always integer multiples of $(1/d!)$.
For convenience we shall call this integer $|\det[1|M]|$ 
the {\em class} of $\sigma$; it is a kind of normalized volume.  
Simplices of class $0$ 
are degenerate.

In the matrix $M$, a choice of any $j+1$ rows corresponds to a
$j$-face $\tau$ 
of $\sigma$.  We shall call these rows the {\em face-rows} of $\tau$;
they represent the vertices of $\tau$.  We call all other rows
of $M$ the {\em non-face-rows} of $\tau$.

In a dual fashion, a choice of any $j$ columns
corresponds to a choice of a $j$-face $F$ of the $d$-cube and all the
$j$-faces parallel to it.  We call these columns the 
{\em cube-face-columns} of $F$;  they correspond to the coordinates
that vary over the face $F$, and outside of these columns, the
coordinates of points on the face $F$ are fixed.

Thus a $j$-face $\tau$ (of a simplex $\sigma$) 
is exterior if and only if there is a
choice of some $j$ columns outside of which the face-rows of $\tau$ are
identical.  Then for exterior faces $\tau$,
we may speak of these columns as the {\em cube-face-columns} of
$\tau$, and all the other columns are the {\em non-cube-face-columns}
of $\tau$.

As an example, the following matrix $M$ represents a simplex $\alpha$ in
the $5$-cube:
\begin{equation}
\label{example-matrix}
\begin{array}{ccccccc}
&\downarrow & \downarrow & \downarrow && &\\
\Rightarrow&0&0&1&1&0&\leftarrow\\
           &1&0&1&1&0&\leftarrow\\
           &0&0&0&1&0&\leftarrow\\
\Rightarrow&0&1&1&0&0&\\
\Rightarrow&0&1&1&1&0&\leftarrow\\
\Rightarrow&0&1&1&1&1&\\
 & &\Uparrow&&\Uparrow&\Uparrow& \\
\end{array}
.
\end{equation}
The rows and columns marked by single-arrows on the right and top
are face-rows and cube-face-columns 
for some exterior $3$-face $\sigma$.  One can verify this
by checking that there are $3$ columns and $3+1$ rows; and 
after deleting the cube-face-columns,
the face-rows look identical.  Another way to say this is in any fixed
non-cube-face-column, the entries in the cube face rows must be the same.
(For instance, in column 4, $m_{41}=m_{42}=m_{43}=m_{45}$.)
Another exterior $3$-face $\tau$ is represented by the face-rows
and cube-face-columns 
marked on the left and
bottom by double-arrows.  It contains an 
exterior $2$-face represented by rows $\{4,5,6\}$ and columns $\{4,5\}$.
The rows $\{1,2\}$ represent an exterior edge because there is a
cube-face-column (i.e., column 1) outside of which the two rows are
identical.  The rows $\{2,3\}$, however, do not represent an exterior
edge.

\section{Constraints on Covers}
A simplicial cover of the $d$-cube induces simplicial covers of each
of its $j$-faces (which are again cubes).
Note also that those covers consist of {\em exterior faces} of the
simplices in the original cover.  Thus, for each dimension $j$, a
natural constraint for a $d$-cube cover is that its
$d$-simplices must have enough $j$-dimensional exterior faces to cover
the $j$-faces of the $d$-cube.

\begin{defn}
Let $F(d,c,d',c')$ count the maximal number of 
dimension $d'$, class $c'$ exterior
faces that any dimension $d$, class $c$ simplex in the $d$-cube can have.
\end{defn}

Given a cover of the $d$-cube, let $x_c$ 
represent the number of simplices of class $c$ in that cover.  We
wish to minimize $\sum x_c$, the total number of simplices in the
cover, subject to some constraints:
\begin{eqnarray}
  \label{eq:constraint}
  \sum_{c=1}^{V(d')}  
        \frac{c}{d'!} F(d,c,d',c) x_c & \geq & 2^{d-d'}{d \choose d'},
        \quad (d' = 1, 2, \ldots, d).
\end{eqnarray}
Recall that $V(d)/d!$ is the volume of the largest $d$-simplex in the
$d$-cube.  Equivalently, $V(d)$ is the class of that largest simplex.

To see how the constraints arise, observe that 
the right side of (\ref{eq:constraint})
counts the total $d'$-volume of the cube-faces of dimension $d'$.
For each dimension $d' \le d$, 
there are $2^{d-d'}{d \choose d'}$ $d'$-faces in the $d$-cube, each
with $d'$-volume 1. 

So there must be enough elements of the cover to cover the $d'$-volume in
each dimension.  Notice that in a cover, 
the $d'$-cube-faces must be covered by exterior facets of the 
exterior $(d'+1)$-faces that cover the $(d'+1)$-cube-faces.  
Thus for $d'$-cube-faces, we need only consider elements of the
cover that arise in a successive chain 
of exterior facets of exterior facets, up
through the top dimension.  Since 
Proposition \ref{class-divides} will show that any exterior facet
of a simplex must have the same class as that simplex,
we need only consider exterior faces that have same class as the
simplices they lie on.
Thus the left side of (\ref{eq:constraint})
counts the maximal volume of the exterior dimension
$d'$-faces that could arise from these elements of the cover, since
$x_c$ is the number
of $d$-simplices of class $c$ in the cover,
$F(d,c,d',c)$ is the maximal number of exterior $d'$-simplices
that a $d$-simplex can have (of the same class $c$), and
$c/d'!$ is the volume of such a simplex.

For these constraints to be helpful, we will need the value of
$F(d,c,d',c')$, or at least, an upper bound, and this allows us to
improve our linear program further later.

\section{Projecting along an Exterior Face}
Let $\alpha$ be a $d$-simplex of class $c$ in the $d$-cube.  
Suppose that $\sigma$ is an exterior $d'$-face of
$\alpha$.  Without loss of generality, we can assume that one of the
vertices of $\sigma$ is at the origin, and that $\sigma$ is contained 
in the $d$-cube-face in which last $d-j$ coordinates are zero.

Consider the matrix representation
$M$ of $\alpha$.  We may assume that the first 
$j+1$ rows of $M$ are the face-rows of
$\sigma$ and the first row is the origin.  
By assumption, the last $d-j$ coordinates of these rows are all zero,
i.e., the cube-face-columns of $\sigma$ are the first $j$ columns of $M$.
Thus $M$ has the following form:
\begin{equation}
\label{matrix-repn}
M = 
\left[
\begin{array}{ccc|ccc}
 0 & \cdots & 0  &   0 & \cdots               & 0 \\ \hline 
   &        &    &     &                      &   \\
   &   A    &    &     &     \mbox{ zeroes }  &   \\
   &        &    &     &                      &   \\ \hline
   &        &    &     &                      &   \\
   &   C    &    &     &  B                   &   \\
   &        &    &     &                      &   \\
\end{array}
\right].
\end{equation}
Here $A$ and $B$ are square submatrices of size 
$j \times j$ and $(d-j) \times (d-j)$ respectively.  The non-degeneracy
of $\sigma$ implies that we can add multiples of the rows containing
$A$ to zero out the submatrix $C$.  This yields a new matrix
$M_\sigma$:
\begin{equation}
\label{matrix-prime}
M_\sigma =
\left[
\begin{array}{ccc|ccc}
 0 & \cdots & 0  &   0 & \cdots               & 0 \\ \hline
   &        &    &     &                      &   \\
   &   A    &    &     & \mbox{ zeroes }  &   \\
   &        &    &     &                      &   \\ \hline
   &        &    &     &                      &   \\
   & \mbox{ zeroes }   &    &     &   B   &   \\
   &        &    &     &                      &   \\
\end{array}
\right].
\end{equation}
Note that $[1|M]$ and $[1|M_\sigma]$ have the same determinant.

Let $\sigma^\perp$ denote the $(d-j)$-simplex spanned by the points
corresponding to the zero
vector and the last $d-j$ rows of $M_\sigma$ 
(containing the submatrix $B$).
Here, $\sigma^\perp$ can be viewed 
as the projection of the simplex $\alpha$ into the orthogonal
complement of $\sigma$. It is clearly an exterior face of the simplex 
$M_\sigma$ of dimension $d-d'$.

Let $\pi_\sigma:[0,1]^d \rightarrow [0,1]^d$ denote the
projection that zeroes out the first $j$ coordinates; 
it collapses the simplex $\sigma$ to the origin, and sends all
other vertices of $\alpha$ to $\sigma^\perp$.  We shall call the map
$\pi_\sigma$ the
{\em projection along the face $\sigma$}.
From this discussion we can draw some immediate consequences:
\begin{prop}
\label{class-divides}
If a $d$-simplex $\alpha$ has an exterior $d'$-face $\sigma$,
then the class of $\sigma$ divides the class of $\alpha$.
In particular, if $\sigma$ has codimension 1, then the class of
$\sigma$ equals the class of $\alpha$. 
\end{prop}

\begin{proof}
Note that the class of $\sigma$ equals $\det(A)$ and $\det(A)$
divides $\det[1|M_\sigma]$, which equals
$\det[1|M]$.  The second statement follows from noting that the
submatrix $B$ of the matrix $M$ that represents $\alpha$ is a $1\times 1$
matrix that must contain a 1 if it is non-degenerate.  Hence the
determinant of $A$ equals the determinant of the whole matrix.
\end{proof}

\begin{prop}
\label{parallel-faces}
Suppose a non-degenerate 
        $d$-simplex $\alpha$ has an exterior $d'$-face $\sigma$
        in the cube-face $F$.  
        Then no cube-face parallel to $F$ contains 
        more than one vertex of $\alpha$.
\end{prop}

\begin{proof}
If a cube-face parallel to $F$ contained two vertices of $\alpha$, 
then two rows of $M$ that intersect the submatrix $B$ 
would agree in the last $(d-j)$ coordinates.
After transforming $M$ to $M_\sigma$ by row operations, these two rows
would be identical.  Hence $\det[1|M_\sigma]=0$, 
contradicting the non-degeneracy of the simplex $\alpha$.
\end{proof}

This immediately implies:

\begin{prop}
\label{no-parallel-exterior-faces}
  A non-degenerate simplex $\alpha$ cannot have exterior $j$-faces in two
  parallel cube-faces, for $j > 0$.  Thus in the matrix
  representation of $\alpha$, a given (non-empty) cube-face-column set
  belongs to at most one exterior face.
\end{prop}

(Note that a class 2 simplex in the 3-cube does have 1-faces in
parallel 2-faces of the 3-cube, but these are not exterior faces by
our definition.)

\begin{prop}
\label{pi-sigma-one-to-one}
The projection $\pi_\sigma$ is one-to-one on the vertices of $\alpha$ that
are not in $\sigma$.
\end{prop}

\begin{proof}
The fiber of the projection $\pi_\sigma$ on any vertex of
$\sigma_\perp$ is a cube-face parallel to $\sigma$; two vertices of
$\alpha$ in one fiber would then contradict Proposition
\ref{parallel-faces}.
\end{proof}

\section{Footprints and Shadows}
Now suppose that, in addition to $\sigma$, 
the simplex $\alpha$ has another exterior $j$-face $\tau$.
In order to count how many such $\tau$ there are, we 
examine what happens to $\tau$ under $\pi_\sigma$, 
the projection along the face $\sigma$.  
Since $\tau$ is a face of $\alpha$, the projection
$\pi_\sigma(\tau)$ is a subset of $\sigma^\perp$.

\begin{defn}
Define the 
{\em shadow} of $\tau$ with respect to $\sigma$ to be
the projection $\pi_\sigma(\tau)$.
Define the {\em footprint} of $\tau$ with respect to $\sigma$ to be the 
intersection $\tau \cap \sigma$.
(Note that the footprint may be empty.)
\end{defn}

\begin{figure}[thpb]
\begin{center}
\includegraphics[height=1.5in]{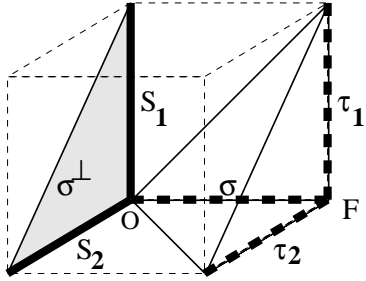}
\end{center}
\caption{A simplex in the cube, with exterior edges denoted by dotted
  lines.  The origin is denoted by $O$. 
  With respect to $\sigma$, the footprint and shadow of $\tau_i$ are
  $F$ and $S_i$, respectively.}
\label{foot-shad}
\end{figure}

Consider the following example.  Let $\alpha$ be the simplex in the cube of 
Figure \ref{foot-shad}.  
The three dotted edges of $\alpha$ are its
exterior $1$-faces: $\sigma$, $\tau_1$, and $\tau_2$.
Suppose that the origin is denoted by $O$.
Then the triangle $\sigma^\perp$ is the projection of the simplex
$\alpha$ into the orthogonal complement of the exterior edge $\sigma$.

With respect to $\sigma$, 
the footprint of $\tau_1$ is the point $F$, 
and the shadow of $\tau_1$ is the edge $S_1$.
Similarly, with respect to $\sigma$, 
the footprint of $\tau_2$ is the point $F$, 
and the shadow of $\tau_2$ is the edge $S_2$.
The exterior face $\sigma$ also has a footprint and shadow with
respect to itself: the footprint is $\sigma$ and the shadow is just
the point at the origin $O$.

The following lemma is crucial in establishing facts about
footprints and shadows.

\begin{lemma}
\label{jk}
Let $M$ be the matrix representation of $\alpha$.  
Suppose that $j$
is the number of face-rows of $M$ that $\sigma$ and $\tau$ have in
common, and $k$ is the number of cube-face-columns that they have in
common.  

(i) If $j>0$ then $j=k+1$.  Else if $j=0$, then $k=0$.

(ii) As long as $k \neq 0$,
the number of {\em non}-face-rows that $\sigma$ and $\tau$ share must
be equal to the number of {\em non}-cube-face-columns they share.
\end{lemma}

As an example, consider the simplex $\alpha$ represented 
by the matrix $M$ in (\ref{example-matrix}) with exterior faces 
$\sigma$ and $\tau$.  In this $\sigma$ and $\tau$ share $j=2$
face-rows (rows 1 and 5) and $k=1$ cube-face-columns (column 2).
Since $k\neq 0$, we expect that the number of non-face-rows they share
(none) should equal the number of non-cube-face-columns they share
(also none).

\begin{proof}
Suppose $\sigma$ has $s+1$ face-rows (and $s$
cube-face-columns) and $\tau$ has $t+1$ face-rows (and $t$
cube-face-columns).

Assume first that $j>0$.
Since the vertices corresponding to common face-rows all lie
in an affine subspace defined by the common cube-face-columns, it must be that
$j\leq k+1$ (else $\alpha$ would be degenerate).
On the other hand, 
the vertices of $\sigma$ and $\tau$ 
together have $s+t+2-j$ face-rows.  
Since $j > 0$, these rows are identical 
outside of $s+t-k$ columns.  Thus $s+t+2-j$ vertices of $\alpha$
lie in an affine subspace of dimension $s+t-k$, so the non-degeneracy of
$\alpha$ means that $s+t+2-j \leq s+t-k+1$, or $j\geq k+1$.  Thus $j=k+1$.

If $j=0$, then $\sigma$ and $\tau$ have no face-rows in common.  
We can assume that face-rows of $\sigma$ are zero
outside its cube-face-columns $\mathcal{C_\sigma}$.  
Let $\mathcal{C_\tau}$ denote the cube-face-columns of $\tau$,
and consider the columns $\mathcal{C_\sigma} \cup \mathcal{C_\tau}$.
Outside of these columns, we zero out the entries in all but one of
the face-rows of $\tau$ by subtracting one of the face-rows from all
the others.  This does not affect the determinant of $M$.
Now the $s+t+2$ face-rows of $\sigma$ and $\tau$ together lie in the affine
subspace determined by $\mathcal{C_\sigma} \cup \mathcal{C_\tau}$,
since they agree outside these columns.  Hence $s+t+1 \leq (s+t-k)+1$,
or $k\leq 0$.  By definition, $k$ cannot be negative, so $k$ must be zero.

For the final assertion, let $j'$ and $k'$ be the number of non-face-rows 
and non-cube-face-columns that $\sigma$ and $\tau$ share in common.
One may check that $j'=d-2d'+j-1$ and $k'=d-2d'+k$.
Then the first assertions show that $j'=k'$ except when $j=0$ 
(and in that case $k=0$ and $k'=j'+1$).
\end{proof}

The lemma has two important corollaries:

\begin{cor}
\label{thm:exterior-intersections}
Let $\sigma$ and $\tau$ be two exterior faces of a non-degenerate
simplex $\alpha$ in the cube.  Then $\sigma \cap \tau$ is an exterior
face, i.e., the footprint of $\tau$ in $\sigma$ is exterior.
\end{cor}

\begin{proof}
If $j=0$ then $\sigma\cap\tau$ is empty, and by definition
exterior.
Otherwise $j>0$, so by Lemma \ref{jk}
$j=k+1$.  This means that the intersection $\sigma\cap\tau$ is of
dimension $j-1=k$.  Since it lies in a cube-face of dimension $k$, it is
an exterior face of $\alpha$.
\end{proof}

\begin{cor}
\label{thm:shadow-exterior}
The shadow of $\tau$ with respect to $\sigma$ is an exterior face of 
$\sigma^\perp$.
\end{cor}

\begin{proof}
If $j=0$ then the shadow $\pi_\sigma(\tau)$ has $t+1$ face-rows, so it
has dimension $t$.
By Lemma \ref{jk}, if $j=0$ then $k=0$, 
so that the $t$ cube-face-columns of $\tau$ are unchanged by the
projection $\pi_\sigma$;
outside of these columns, the face-rows of $\pi_\sigma(\tau)$ are
identical.  So the shadow of $\tau$ is contained in
a cube-face of dimension $t$, as desired.

If $j>0$, then the shadow of $\tau$ has $t+1-j$ non-zero face-rows and
one zero face-row (since vertices of $\sigma$ map to $0$ by
$\pi_\sigma$) and hence the shadow of $\tau$ has dimension $t+1-j=t-k$.
But Lemma \ref{jk} certifies that $j=k+1$, so 
then the shadow of $\tau$ is contained in a cube-face of dimension
$t-k$ as well.
\end{proof}

For our purposes, 
the most important property of footprints and shadows is the following:
\begin{prop}
\label{unique}
Given a simplex $\alpha$ in the cube, fix an exterior $d'$-face $\sigma$.
Then any exterior $d'$-face $\tau$ of $\alpha$ has a unique
footprint-shadow pair with respect to $\sigma$.
\end{prop}

(Recall that this footprint may be empty.)

\begin{proof}
By Proposition \ref{pi-sigma-one-to-one}, 
under $\pi_\sigma$,
every vertex of $\sigma^\perp$
(except for the origin) has a unique pre-image in
$\alpha$.  Thus the vertices of $\tau$ not in $\sigma$ can be
determined from the shadow of $\tau$.  Clearly, the vertices of $\tau$ in
$\sigma$ can be determined from the footprint of $\tau$.
\end{proof}

\section{A recursion for $F(d,c,d',c')$}
Given a $d$-simplex $\alpha$ of class $c$ in the $d$-cube, 
we wish to count how
many exterior faces it may have of dimension $d'$, class $c'$.
Call such a face an exterior $(d',c')$-face of $\alpha$.  
If it cannot have any, then $F(d,c,d',c')=0$.  Otherwise, it has at
least one, call it $\sigma$.

Proposition \ref{unique} shows that once the exterior face $\sigma$ is
fixed, any other exterior $d'$-simplex $\tau$ determines
a unique footprint-shadow pair.  The footprint and shadow are exterior faces 
(of $\sigma$ and $\sigma^\perp$ respectively),
their dimensions must add up to $d'$,
and the products of their classes must be $c'$.  
(The previous sentence holds for empty footprints if empty sets are
considered to have dimension 0 and class 1.)
So the maximal number of exterior $(d',c')$-faces 
must be bounded above by the number of footprint-shadow pairs, summed over
the dimension and class of the footprint.

\begin{thm}
  \label{thm:recurrence}
If $c>V(d)$, $c'>V(d')$, $d'>d$, or $c' \not | c$, then
$F(d,c, d',c')=0$.  Otherwise,  $F$ obeys the recursion
 \begin{equation}
   \label{eq:recurrence}
   F(d,c,d',c') \leq \sum_{\delta=0}^{d'}\sum_{\gamma=1}^{c'}
        \, F(d',c',\delta,\gamma)
        \cdot F(d-d',c/c',d'-\delta,c'/\gamma),
 \end{equation}
 where $F(d,c,0,1)$ is taken to be $1$.
\end{thm}

\begin{proof}
If $c>V(d)$, $c'>V(d')$, $d'>d$, or $c' \not | c$, then the definition
of $V(d)$ and Proposition \ref{class-divides} 
rule out the possibility of any exterior $(d',c')$-faces.  Otherwise, if there
is an exterior $(d',c')$-face $\sigma$, then
for any exterior $(d',c')$-face $\tau$ consider the
footprint and shadow of $\tau$ with respect to $\sigma$.

Let $\delta$ be the dimension of the footprint, and let $\gamma$ be
its class.  Then $F(d',c',\delta,\gamma)$ counts the maximal number of ways in
which the footprint of $\tau$ could be an exterior face of $\sigma$.
Similarly, because of Proposition \ref{class-divides},
$F(d-d',c/c',d'-\delta,c'/\gamma)$ counts the maximal number of
ways that the shadow of $\tau$ can be an exterior face of
$\sigma^\perp$.  Since footprint-shadow pairs are unique 
(Proposition \ref{unique}),  summing over all possible footprint
dimensions, the number of such pairs that $\tau$ could
assume is given by the right hand side of Equation
(\ref{eq:recurrence}).

The case where the footprint is empty merits special caution.
In this case, we choose the convention that the
footprint be $0$-dimensional (not $(-1)$-dimensional, as one might
expect).  There are two reasons we do this.   First, there is no overlap with
the footprint-shadow pairs that arise from $\delta=0$, because if
$\tau$ has empty footprint and the footprint of $\tau'$ is a vertex
$v$, then
the shadows of $\tau$ and $\tau'$ must be different
($\pi_\sigma(v)=0$, so one shadow contains $0$ and the other does not).
Secondly, in both cases the dimension of the shadow must be $d'$.

This convention simplifies the sum (\ref{eq:recurrence}), removing the
need for an extra term for $\delta=-1$.
\end{proof}

The recursion yields some closed form expressions
for $F(d,c,d',c')$ in some special cases.
Let $\Delta(c) = \min \{d:V(d) \ge c \}$.  Hence the
smallest dimension in which a class-$c$ simplex appears is at least
$\Delta(c)$, with equality if $c=V(d)$.
Table~\ref{delta-table} shows some known values of $\Delta(c)$.
\begin{table}[htbp]
  \caption[Some values of $\Delta(c)$]
        {Some values of $\Delta(c)$.  Compare Table \ref{vd-table}.}  
  \label{delta-table}
  \begin{tabular}{|r|lllllllllllll|}
    \hline
    $c$ & 1&2&3&4&5&9&32&56&144&320&1458&3645&9477\\
    $\Delta(c)$ &0&3&4&5&5&6&7&8&9&10&11&12&13\\
    \hline
  \end{tabular}
\end{table}

\begin{thm}
\label{c-prime=c}
If $c'=c$, then 
\begin{equation}
\label{good-bound}
F(d,c,d',c) \leq {d-\Delta(c) \choose d'-\Delta(c)}.
\end{equation}
\end{thm}
%
\begin{proof}
We prove the theorem by induction on $k=d-\Delta(c)$.  
If $k=0$, then $d=\Delta(c)$.  Hence for any $d'<\Delta(c)$, 
both sides of (\ref{good-bound}) are 0,
because $\Delta(c)$ is a lower bound for the 
dimensions in which a class-$c$ simplex appears.
If $d'=\Delta(c)$, then the left side of (\ref{good-bound}) is at most
1, and the right side is 1, confirming the desired inequality.

Now assume the theorem holds whenever $d-\Delta(c) \leq k-1$.  
We shall show it holds for $d-\Delta(c)=k$.

If $d'=d$, then it is easy to check that both sides of 
(\ref{good-bound}) are 1, as desired, and 
if $d'<\Delta(c)$, then both sides are 0, as desired.
The remaining cases for $d'$ are if $\Delta(c)\leq d'< d$.
Using the recurrence of Theorem \ref{thm:recurrence}, 
the only nonzero terms in the sum have $\gamma=c$ since $c'=c$,
and $\Delta(c) \leq \delta \leq d'$, so we have:
\begin{equation}
\label{reduced-recurrence}
F(d,c,d',c)  \leq  \sum_{\delta=\Delta(c)}^{d'} 
                 F(d',c,\delta,c) \cdot F(d-d',1,d'-\delta,1).
\end{equation}
Then if $d'=\Delta(c)$, this inequality becomes
$$F(d,c,d',c) \leq F(\Delta(c),c,\Delta(c),c) \cdot
F(d-\Delta(c),1,0,1).$$
By definition, $F(d-\Delta(c),1,0,1)=1$ and
$F(\Delta(c),c,\Delta(c),c)$ is $1$ if a class $c$ simplex can be realized
in dimension $\Delta(c)$ and it is $0$ otherwise.
In any case, their product is less than or equal to 
${d-\Delta(c) \choose \Delta(c)-\Delta(c)}=1$, as desired.
Otherwise, $\Delta(c)< d'< d$, so the inductive hypothesis
can be used on the right side of (\ref{reduced-recurrence}), since
$d'-\Delta(c)$ and $d-d'-\Delta(1)$ are strictly less than $d-\Delta(c)$.
Hence:
$$
  F(d,c,d',c) 
       \leq \sum_{\delta=\Delta(c)}^{d'} 
      {{d'-\Delta(c)} \choose {\delta-\Delta(c)}} {{d-d'} \choose {d'-\delta}} 
      = {{d-\Delta(c)} \choose {d'-\Delta(c)}},
$$
as desired.
\end{proof}

In particular, when $c'=c=1$, the bound in Theorem \ref{c-prime=c} becomes 
${ d \choose d'}$ and it is achieved by the corner simplices of the $d$-cube,
i.e., simplices spanned by one vertex and all its
nearest neighbors.  In fact, as we show in the following theorem,
if $1<d'<d$, corner simplices are the only
simplices for which the bound on $F(d,1,d',1)$ is sharp.
\begin{thm}
\label{non-corner}
If $1<d'<d$ and $c=1$, then the bound of Theorem \ref{c-prime=c} 
on $F(d,1,d',1)$
is sharp for corner simplices.  
In fact, for $1<d'<d$, a simplex has strictly more than 
$\frac{d-1}{d} {d \choose d'}$ exterior
$d'$-faces if and only if it is a corner simplex.  Hence for
non-corner simplices, the number of exterior $d'$-faces is
bounded above by $\lfloor \frac{d-1}{d} {d \choose d'} \rfloor$.
\end{thm}

Note that the ``only if'' implication does not hold if $d'=1$ or if $d'=d$;
for example if $d=3$ and $d'=1$, then a non-corner simplex can have 
as many exterior edges as a corner (three), and if $d=d'$, then any
simplex is its own exterior $d'$-face.  Note also that the bound 
$\lfloor \frac{d-1}{d} {d \choose d'} \rfloor$ in
the final assertion holds even for non-corner simplices of class
$c>1$, although in such cases, 
the bound of Theorem \ref{c-prime=c} is better.

\begin{proof}
Let $\alpha$ denote a corner simplex. Without loss of generality its
matrix representation $M$ can be written as the identity $d \times d$
matrix augmented with a row of zeroes.  Any choice of $d'$ columns of
the $d$ columns of $M$ specifies a set of cube-face-columns, and one
may easily locate $d'+1$ rows of $M$ outside of which those 
cube-face-columns are identical.  These are face-rows for some
exterior $d'$-face of $\alpha$.  Thus there is one such face (and by
Proposition \ref{no-parallel-exterior-faces}, at most one) for each
of the ${ d \choose d'}$ choices of cube-face-column sets.

For the converse, let $M$ be a matrix representation of $\alpha$, and
let $\mathcal{C}_j$ denote the set of $d-d'$ consecutive columns of $M$
that start at the $j$-th column (and wrap around if necessary).  
For instance, if $M$ has $d=8$ columns and $d'=3$, then 
$\mathcal{C}_3$ would be the column set $\{3,4,5,6,7\}$, and 
$\mathcal{C}_6$ would be column set $\{6,7,8,1,2\}$.  Let $\mathcal{R}_i$ be
defined similarly for rows of $M$ {\em excluding the $(d+1)$-st row},
e.g., in the previous example, even though $M$ has $d+1=9$ rows, 
$\mathcal{R}_6$ would be row set $\{6,7,8,1,2\}$.

We first show that if the number of exterior $d'$-faces of $\alpha$
is greater than $\frac{d-1}{d} {d \choose d'}$, then $\alpha$ has a 
{\em $d$-cycle of exterior $d'$-faces}, in other words,
there must be an ordering of the columns of $M$ such that 
each of $\mathcal{C}_1, \mathcal{C}_2, ... , \mathcal{C}_d$ 
represents non-cube-face-columns for some exterior $d'$-face.

Let $d^*=d-d'$.
There are $d!$ possible orderings of the columns of $M$, and we claim
that at least 
one of them will remain as a possible ordering for a $d$-cycle if at
most $1/d$ of the possible ${d \choose d'}$ column subsets of size
$d^*$ are {\em not} non-cube-face-columns for some exterior face.  This is true
because if some subset $A$ of $d^*=d-d'$ columns are {\em not}
non-cube-face-columns for some exterior face, then this restriction
rules out $(d-d')! d'! d$ of the $d!$ orderings of the columns of $M$
(since there are $(d-d')!$ ways to order $A$, $d'!$ to order the
columns not in $A$, and $d$ ways to place the columns of $A$ adjacent
to each other).  Thus if we rule out strictly fewer than 
$d!/ (d-d')! d'! d = {d \choose d'}/d$ exterior faces, there will still
remain an ordering that could occur as an ordering of the columns of
$M$, and in such an ordering, every $\mathcal{C}_i$ is a set of
non-cube-face-columns for some exterior face.

Now we show that if $\mathcal{C}_1,...,\mathcal{C}_d$ all represent
non-cube-face-column sets for exterior $d'$-faces, 
then $\alpha$ must be corner simplex. 

Consider $\mathcal{C}_1$, the first $d^*$ columns of
$M$; by assumption, these are 
non-cube-face-columns for some exterior $d'$-simplex $\sigma_1$.  
We can assume (by re-ordering rows if needed) that its non-face-rows
are $\mathcal{R}_1$, the first $d^*$ rows of $M$.
The intersection of the non-face-rows and
non-cube-face-columns forms a $d^* \times d^*$ block $B_1$ in the matrix $M$.
In any single non-cube-face-column, all entries not in $B_1$ must be
identical.  Note that, by symmetry, toggling all the column elements 
in a column (exchanging 1's and 0's) does not change the isomorphism
class of the simplex, and we can use this operation, if needed, in columns 
$1$ through $d^*$ so that entries not in $B_1$ are all 0.

Similarly, consider columns in $\mathcal{C}_2$
(columns 2 through $d^*+1$);  these are non-cube-face-columns
corresponding to an
exterior $d'$-simplex $\sigma_2$.  Because $d'>1$, $\sigma_1$ and $\sigma_2$ 
share a common cube-face-column (e.g., column $d$), so
Lemma \ref{jk}(ii) applies: noting that that 
$\sigma_1$ and $\sigma_2$ share $d^*-1$ non-cube-face-columns, 
they must also share $d^*-1$ non-face-rows.  
We can thus reorder the rows of $M$ so that
$\sigma_2$ occupies rows in $\mathcal{R}_2$ (rows 2 through $d^*+1$).  
Thus the corresponding $d^* \times d^*$ block $B_2$ of $\sigma_2$ 
intersects the block $B_1$ of $\sigma_1$ 
in a block along the diagonal, and as before,
we may toggle column elements in non-cube-face-columns of $\sigma_2$ 
so that all entries in those columns but outside $B_2$ are 0.

In the same way, by cycling through the columns, we can infer from
Lemma \ref{jk}(ii) that each pair of neighboring faces in the 
$d$-cycle of exterior $d'$-faces must share $d^*-1$ non-face-rows
(because they share $d^*-1$ non-cube-face-columns and $d'>1$).  In
particular, since their are $d$ such faces in the cycle, 
there are only $d$ rows of $M$ that are used as non-face-rows of such faces.  
Thus there is one row of $M$ that is not a non-face-row of any face,
in other words, it is a face-row of {\em every} face in the cycle.

We can reorder the rows of $M$ in such a way that this common face-row is 
in the $d+1$-st row of $M$, and the exterior face whose 
non-cube-face-columns are $\mathcal{C}_k$ has non-face-rows 
$\mathcal{R}_k$, for all $1\leq k \leq d$.
Then, by toggling columns if needed, each of the corresponding 
$d^* \times d^*$ blocks along the diagonal now has zeroes in its
columns outside it.  This forces row $d+1$ to be all zeroes, and all
non-diagonal entries are also all zero.  Then all diagonal entries
must be 1 because $\alpha$ is non-degenerate.  Hence $\alpha$
must be isomorphic to a corner.
\end{proof}

\section{Improving the Linear Program}

We can now improve the linear program in (\ref{eq:constraint}) by
noting that Theorem \ref{c-prime=c} shows that all the coefficients
$F(d,c,d',c)$ of $x_c$ are equal for values of $c$ between any two
values of $V$, i.e., $V(k-1)<c\leq V(k)$.  Since the objective
function also equally weights the variables $x_c$, there is an optimal solution
to the linear program (\ref{eq:constraint}) whose support lives
entirely on the variables $x_{V(k)}$ for $k=1,...d'$, e.g., $x_1, x_2,
x_3, x_5, x_9$, etc.
Setting $y_k=x_{V(k)}$ for $k\geq 2$, the program becomes:
\begin{eqnarray}
  \min \sum_{k=2}^{d} y_k  \mbox{   subject to } & & \nonumber\\
  \sum_{k=2}^{d'}  
        \frac{V(k)}{d'!} F(d,V(k),d',V(k))\ y_k & \geq & 2^{d-d'}{d \choose d'}, 
        \quad d' = 1, 2, \ldots, d.
\label{general-program}
\end{eqnarray}

We can improve this linear program further by
considering the presence of corner simplices separately from other possible
exterior class-1 $d'$-simplices.  In the program above, replace $y_2$
(which counts the total number of class-1 simplices) by $y_1 + y_2$, where 
$y_1$ denotes the number of corner simplices in a cover, and $y_2$ now
denotes the number of class-1 non-corners.  Applying 
Theorems \ref{c-prime=c} and \ref{non-corner},
and the observation that $y_1$ must be bounded above by the
total number of vertices of the $d$-cube (since there is at most one
corner simplex per vertex), we obtain the program:
\begin{eqnarray}
  \min \sum_{k=1}^{d} y_k  \mbox{   subject to } \nonumber\\
                                     y_1 & \leq & 2^d, \nonumber\\
        \sum_{k=1}^{d'}  
        \frac{V(k)}{d'!} 
\lfloor 
\epsilon(k,d) 
{d- \Delta(V(k)) \choose d'- \Delta(V(k))}
\rfloor
y_k & \geq & 2^{d-d'}{d \choose d'},
        \quad d' = 1, 2, \ldots, d,
\label{reduced-program}
\end{eqnarray}
where $\epsilon(k,d)=1$ unless both $k=2$ and $1<d'<d$, 
in which case $\epsilon(k,d)= \frac{d-1}{d}$.
Note that $\Delta(V(k))$ is almost always $k$, unless $k\leq 2$, and
the floor function only comes into play when $\epsilon(k,d) \neq 1$.

\section{Conclusion}

We solved the linear program (\ref{reduced-program}) 
using the package {\tt lp\_solve},
and rounded non-integer optimal values up to integers.
This yields the bounds in Table \ref{table-bounds}.

It should be possible to improve our bounds for covers by using
more information about the kinds of simplices that can occur in the
cube and how they fit together, similar to what Hughes \cite{hughes} and 
Hughes-Anderson \cite{hughes-anderson} did to study 
$D^v,T^v$-minimal triangulations.
These programs rely on enumerating specific features
of configuration classes of simplices in the cube, or on enumerating
the classes themselves. In the latter case, the resulting linear program 
has one variable for each isomorphism class of configurations, which
for $d=6$ already involves thousands of variables.  
The Hughes approach \cite{hughes}
does not work beyond $d>11$ and the 
Hughes-Anderson program \cite{hughes-anderson} 
becomes intractable for $d>7$.

By contrast, the number of variables in our program for $C(I^d)$
is just $d$.  Also, our results for $C(I^d)$ using a comparatively 
small program compares quite well to their lower bounds for $T^v(I^d)$
and does not require knowledge of the configuration classes that occur 
in specific dimensions (other than the existence of 
corners, which occur in every dimension).  So our program can be
solved for much larger values of $d$ than programs 
that use extra information about configuration classes in specific
dimensions.

Our bounds dominate Smith's bound up through at least $d=12$, and 
computational evidence (limited by roundoff errors)
show that our bounds dominate Smith's bound (\ref{smith-bound}) for covers and 
general triangulations in dimensions up through at least $d=27$,
using specific values of $V(d)$ from Table \ref{vd-table} and the asymptotic
bound for $V(d)$ from (\ref{vd-upperbound}).  (Bounds for higher $d$
could not be computed due to solver overflow errors.)  
However, we do not believe that our bounds will exceed Smith's bound 
asymptotically.

Our bounds do give some new observations in dimension 4.  
It was known \cite{cottle} that the minimal triangulation of the 4-cube
has size 16, for triangulations using only the vertices of the cube.
Our bounds and Theorem \ref{covers-bound} 
then give a stronger result:

\begin{thm}
\label{4-cube-cover}
The minimal simplicial cover of the 4-cube has size 16, and is achieved by a
triangulation.  Consequently, adding vertices to the 4-cube will not 
reduce the size of the minimal triangulation.
\end{thm}

In higher dimensions,
it may be the case that the size of a minimal simplicial
cover of a $d$-cube is strictly smaller than size of the minimal
triangulation.  
(Such a cover may have more symmetry than the minimal triangulation.)
The intuition here is that the largest simplices in
high-dimensional cubes tend to overlap, but more coverage might be
possible by using such 
large simplices than by using small non-overlapping ones.

\begin{problem}
Find a $d$-cube for which the
minimal cover is strictly smaller than the minimal $T^v$-triangulation.
\end{problem}

If this is true, it would suggest that there may be
a $d$-cube for which adding extra vertices helps to make the
triangulation smaller.

\begin{problem}
Find a $d$-cube for which the
minimal $T^v$-triangulation is reduced in size by adding extra vertices.
\end{problem}

One may also ask a similar question about dissections, and note from
Table \ref{table-bounds} that
even in dimension 4 a question remains:

\begin{problem}
Find a $d$-cube for which the minimal dissection is strictly smaller
than the minimal triangulation.  
Does the minimal dissection of the $4$-cube have
15 simplices or 16 simplices?
\end{problem}

We have been pleasantly surprised by the richness of the geometry of cubes.

\subsubsection*{Acknowledgements}
The authors are grateful to G\"unter Ziegler and 
Francisco Santos for helpful conversations.

\bibliographystyle{abbrv}
\bibliography{cube}

\end{document}